\documentclass{article}
\usepackage{graphicx} % Required for inserting images
\usepackage{abstract}
\usepackage{appendix}
\usepackage{hyperref} 
\hypersetup{
hidelinks,
colorlinks=true,
linkcolor=black,
urlcolor=black,
citecolor=black
}
\usepackage{amsmath}
\usepackage{amssymb}
\usepackage[numbers,sort&compress]{natbib}
\usepackage{enumerate}
\usepackage[a4paper, total={6in, 8in}]{geometry}
\numberwithin{equation}{section}
\usepackage{booktabs}

\usepackage{cases}
\usepackage{indentfirst}
\usepackage{algorithmic}
\usepackage[linesnumbered, ruled]{algorithm2e}
\usepackage{ntheorem}
\usepackage{mathtools}
\usepackage{multirow}
\usepackage{slashed}
\usepackage{subfig}

\SetKwRepeat{Do}{do}{while}
\newtheorem{thm}{\textit{Theorem}}[section]
\newtheorem{pro}{\textit{Proposition}}[section]
\newtheorem*{myproof}{\textit{Proof}}
\newtheorem{mydef}{\textit{Definition}}[section]
\newtheorem{cor}{\textit{Corollary}}[section]
\newtheorem{lem}{\textit{Lemma}}[section]
\newtheorem{ass}{\textit{Assumption}}[section]
\newtheorem{rem}{\textit{Remark}}[section]
\RequirePackage{amssymb}
\title{\textbf{Nonconvex ADMM for Rank-Constrained Matrix Sensing Problem}}

\author{Zekun Liu \\School of Mathematical Sciences, Shanghai Jiao Tong University}

\date{}

\begin{document}

\maketitle

\begin{abstract} 
Low-rank matrix approximation (LRMA) has been arisen in many applications, such as dynamic MRI, recommendation system and so on. The alternating direction method of multipliers (ADMM) has been designed for the nuclear norm regularized least squares problem and shows a good performance. However, due to the lack of guarantees for the convergence, there are few ADMM algorithms designed directly for the rank-constrained matrix sensing problem (RCMS). Therefore, in this paper, we propose an ADMM-based algorithm for the RCMS. Based on the
Kurdyka-Łojasiewicz (KL) property, we prove that the proposed algorithm globally converges. And we discuss a specific case: the rank-constrained matrix completion problem (RCMC). Numerical experiments show that specialized for the matrix completion, the proposed algorithm performs better when the sampling rate is really low in noisy case, which is the key for the matrix completion.    \\
\textbf{keywords} Rank-constrained matrix snesing, ADMM, global convergence, matrix completion
\end{abstract}

\section{Introduction}\label{sec:1}

Low-rank matrices arise widely in the field of computer science and applied mathematics, such as the Netflix problem \cite{netflix}, quantum state tomography \cite{quantum}, sensor localization \cite{sensor1,sensor2} and natural language processing \cite{nlp1,nlp2}. In all of these applications, the data size can be extremely large, hence it is expensive and impossible to fully sample the entire data, which leads to the incomplete observations. A natural question is how can we recover the original data from these incomplete observations. In general, such recovery is not guaranteed to be possible. However, with the prior information that the matrix is low-rank \cite{lrcc}, it is possible to recover the original matrix from incomplete observations in an efficient way. One of the most common and widely used model in matrix recovery is the low-rank matrix approximation \cite{lrma1,lrma2}.

For the general low-rank matrix approximation problem (LRMA), we will only observe $b=\mathcal{A}(X)+e$ where $e$ consists of independent white Gaussian noise and $\mathcal{A}:\mathbb{R}^{m\times n}\to \mathbb{R}^{d}$ is the linear measurement operator which transforms $X$ to $(\left \langle X,A_1  \right \rangle,\left \langle X,A_2  \right \rangle,\cdots,\left \langle X,A_d  \right \rangle)^{\top}$ with $A_1,A_2,\cdots,A_d \in \mathbb{R}^{m\times n}$ given and $\left \langle \cdot,\cdot  \right \rangle$ being matrix inner product. In this case, the formulation of the low-rank approximation problem is the following rank-constrained matrix sensing problem (RCMS):
\begin{equation}\label{equ:1.1}
\begin{aligned}
        \min_{X\in \mathbb{R}^{m\times n}} &\left \| \mathcal{A}(X)-b \right \|_{2}^{2} \\
        s.t. \enspace &rank(X)\le r,
\end{aligned}
\end{equation}
where $r<\min\left \{ m,n \right \}$ denotes the upper tight estimation of the rank of the target matrix.

As stated in \cite{ADMiRA}, the formulation (\ref{equ:1.1}) can handle both noiseless and noisy case in a single form. And it can also work for the case that the original matrix is not exactly low-rank but can be approximated accurately by a low-rank matrix. Besides, comparing to the fixed-rank constraint, low-rank constraint can avoid the uncertainty about the existence of the solutions and the complexity of the algorithm convergence analysis 
which are both arisen from
the closure of the fixed-rank constraint \cite{fixrank}. 

\subsection{Related work}
Due to the combinatorial characteristic of rank, problem (\ref{equ:1.1}) belongs to the NP-hard combinatorial optimization problem, which is difficult to analyse in theory and design algorithms \cite{fixrank}. A widely used convex relaxation to the rank minimization problem is the nulear norm (NN) given by \cite{nuclear}. Based on this convex relaxation, singular value thresholding (SVT) \cite{SVT} is proposed to solve the specific low-rank matrix completion problem (LRMC) where the operator $\mathcal{A}$ only takes a subset of the original matrix \cite{lrmc}. However, the efficiency of SVT is only guaranteed in the noiseless case with linear equality constraints. Inspired by the similarity and connection between rank minimization for matrices and $\ell_0$-norm minimization for vectors \cite{relation}, atomic decomposition for minimum rank approximation (ADMiRA) \cite{ADMiRA} is proposed, which is a generalization of the compressive sampling matching pursuit (CoSaMP) \cite{CoSaMP} to matrices. And it has strong theoretical guarantees. But for RCMC, because the linear operator does not satisfies the R-RIP \cite{lrcc}, ADMiRA has no theoretical guarantees. Iterative hard thresholding (IHT) \cite{IHT1,IHT2} can also work for the rank-constrained minimization. And it converges really fast when the rank of the original matrix is very low. Specialized for LRMC, some manifold optimization algorithms such as \cite{mani1,mani2} are efficient and theoretically sound. 

As a powerful and widely used first-order algorithm, the alternating direction method of multipliers (ADMM) \cite{OriADMM1,OriADMM2} also shows a good performance for the LRMA. For the convex nuclear norm regularized least squares problem, ADMM not only performs excellently in practice, but also has a theoretical guarantee for convergence. Considering the specific matrix completion problem, \cite{FADMM} proposes the truncated nuclear norm (TNN) as a better approximation of the rank than the nuclear norm, and uses the ADMM with adaptive penalty to solve their relaxation problem. Moreover, \cite{ADMM1} considers the joint weighted and TNN to approximate the rank with higher accuracy, and applies the ADMM on the relaxation problem for the matrix completion-assisted mmWave MIMO channel estimation. However, both of them use a rank relaxation and only consider the specific matrix completion problem. In the application of model identification, \cite{selfad} directly applies the ADMM to the rank-constrained model identification instead of using a relaxation of rank, and proposes a self-adaptive penalty method. However, due to the rank-constraint, they do not talk about the convergence of the nonconvex ADMM.
\cite{NADMM} uses the matrix factorization to deal with the rank-constraint, and proves the three-block ADMM for matrix completion problem converges to the KKT point. Similarly, \cite{MFADMM} also handles the rank-constraint with matrix factorization, and they apply the four-block ADMM to the general semidefinite programming (SDP) problem instead of the specific matrix completion. Besides, they prove the algorithm converges to the KKT point of their problem. However, since the observations are limited on the equality constraint, the problem formulation in \cite{MFADMM} cannot handle the noisy case as well as the observations being at the objective function like (\ref{equ:1.1}).

\subsection{Contributions}

Since the rank-constraint of matrices is an extension of $\ell_{0}$-norm constraint of vectors, we desire to generalize our previous research on the multiple measurement vector (MMV) problem with the $\ell_{2,0}$-norm \cite{MMVADMM} to the RCMS. 
Unlike most other work \cite{SVT,FADMM,ADMM1,NADMM} focusing on the specific matrix completion problem, 
in this paper, we directly study the general matrix sensing problem, where the measurement operator $\mathcal{A}$ is linear. Moreover, instead of using a relaxation of the rank such as NN \cite{nuclear}, TNN \cite{FADMM} and so on, or handling the rank-constraint with the matrix factorization like \cite{NADMM,MFADMM}, we directly focus on the original rank-constraint problem (\ref{equ:1.1}), whose advantages have been stated in the beginning.

Through applying ADMM to the reformulated problem, we obtain the algorithm called RCMS-ADMM for solving the rank-constrained matrix sensing. As a specific application of the matrix sensing, we also talk about some details of RCMS-ADMM for the significant matrix completion problem. And cause the rank-constraint is nonconvex, we devote length to the description of RCMS-ADMM being globally convergent instead of omitting the theoretical convergence analysis. Moreover, to our knowledge, we first figure out the relation between the global optimal solution of problem (\ref{equ:1.1}) and the Lagrangian multiplier. Besides, for the specific matrix completion problem, compared with the similar algorithms, experiments show that the proposed algorithm has an advantage in reconstruction accuracy when the sampling rate is really low, which is the key of matrix completion.

\section{Preliminaries}\label{sec:2}

In this section, we list some notations and definitions which are used for further analysis.

Notations: $rank(X)$, $X^{\top}$ and $\left \| X \right \|_F$ represent the rank, the transpose and the Frobenius-norm of the matrix $X$ respectively. $\left \langle \cdot,\cdot \right \rangle$ denotes the inner product of two matrices of equal size. $I_k$ represents the $k\times k$ identity matrix. While $I$ denotes the identity reflection. $\mathbf{1}$ represents the matrix of all ones. $vec(X)$ denotes the vector given by concatenating each column of the matrix $X$ in order. $x_m$ represents the $m$-th element of the vector $x$. And $a_{ij}$ denotes the $(i,j)$-element of the matrix $A$ . $A\odot B$ denotes the element-wise product of two matrices of equal size. While $A\oslash B$ denotes the element-wise division of two matrices of equal size. $A\succeq B$ means that the matrix $A-B\in \mathbb{S}^n$ is semi-positive definite. For two composable operators $\mathcal{A}$ and $\mathcal{B}$, $\mathcal{A}\mathcal{B}$ represents their composition $\mathcal{A}\circ \mathcal{B}$.

\begin{mydef}[\cite{view1}]\label{def:2.1}
    For the linear operator 
    \begin{align*}
    \mathcal{A}:&\mathbb{R}^{m\times n}\to \mathbb{R}^d  \\
    &X\mapsto (\left \langle A_1,X \right \rangle,\left \langle A_2,X \right \rangle,\cdots,\left \langle A_d,X \right \rangle)^{\top},
    \end{align*}
    the adjoint of $\mathcal{A}$ is defined as $\mathcal{A}^{\ast}(w)=\sum_{i=1}^d w_i A_i$, where $w\in \mathbb{R}^d$.
\end{mydef}

\begin{mydef}[\cite{cop}]\label{def:2.3}
    For the generalized real function $f:\mathbb{R}^{m\times n}\to \mathbb{R}\cup \left \{ \pm \infty \right \}$.
    \begin{enumerate}[(i)]
    \item Given a nonempty set $\mathcal{X}$, call $f$ proper to $\mathcal{X}$ if $\exists x\in \mathcal{X}$ such that $f(x)<+\infty$, and $\forall x\in \mathcal{X}, f(x)>-\infty$.
    \item $f$ is lower semicontinuous if $\forall x\in \mathbb{R}^{m\times n}, \liminf_{y\to x}f(y)\ge f(x)$.
    \item $f$ is closed if its epigraph 
    \begin{displaymath}
        epif=\left \{ (x,t)\in \mathbb{R}^{m\times n}\times \mathbb{R}|f(x)\le t\right \}
    \end{displaymath}
    is closed.
    \item $f$ is Gradient-$L$-Lipschitz continuous if $\exists L>0$, for $\forall x,y\in \mathbb{R}^{m\times n}$, $\left \| \nabla f(x)-\nabla f(y) \right \| \le L\left \| x-y \right \|$.
\end{enumerate}
\end{mydef}

\begin{pro}[\cite{cop}]\label{pro:2.1}
    For the generalized real function $f$, $f$ is lower semicontinuous iff $f$ is closed.
\end{pro}

\begin{mydef}[\cite{semi}]\label{def:2.4}
    $S$ is a semi-algebraic set in $\mathbb{R}^{m\times n}$ if there exist polynomials $g_{ij},h_{ij}:\mathbb{R}^{m\times n}\to \mathbb{R}$ with $1\le i\le p,1\le j\le q$, such that 
    \begin{displaymath}
        S=\bigcup_{i=1}^{p}\bigcap_{j=1}^{q}\left \{x\in  \mathbb{R}^{m\times n}:g_{ij}(x)=0,h_{ij}(x)>0 \right \}.
    \end{displaymath}
    And for the proper function $f$, it is semi-algebraic if its graph
    \begin{displaymath}
        graphf=\left \{(x,t)\in \mathbb{R}^{m\times n}\times\mathbb{R}:f(x)=t \right \}
    \end{displaymath}
    is a semi-algebraic set in $\mathbb{R}^{m\times n}\times\mathbb{R}$.
\end{mydef}

\begin{pro}[\cite{semi}]\label{pro:2.2}
    Semi-algebra has the following properties.
    \begin{enumerate}[(i)]
        \item Semi-algebra is stable under finite unions.
        \item Real polynomials are all semi-algebra.
        \item The indicator functions of semi-algebraic sets are also semi-algebra.
        \item Suppose $A$ is a semi-algebra in $\mathbb{R}^{m\times n}\times\mathbb{R}$, then its projection $\pi (A)$ where $\pi:\mathbb{R}^{m\times n}\times\mathbb{R}\to \mathbb{R}^{m\times n}$ is also a semi-algebra.
    \end{enumerate}
\end{pro}

\begin{pro}[\cite{semiKL}]\label{pro:2.3}
If $f$ is proper lower semicontinuous and semi-algebraic, then it is also a KL function.
\end{pro}

As for the definition of subdifferential, normal cone, Karush-Kuhn-Tucker (KKT) conditions and Kurdyka-Łojasiewicz (KL) property, we recommend the interested readers to refer to \cite{cop,semi,semiKL}.

\section{Problem formulation}\label{sec:3}

In order to apply ADMM to problem (\ref{equ:1.1}), we rewrite the rank-constrained problem (\ref{equ:1.1}) to make it a two-block optimization problem with a linear equation constraint. 

First, introduce the indicator function 
\begin{displaymath}
    \delta_{\mathcal{C}}(X)=\begin{cases}
 0, & \text{ \textit{if} } X\in \mathcal{C} \\
 +\infty, & \text{ \textit{if} } X\notin \mathcal{C}
\end{cases}
\end{displaymath}
of the rank-constrained set $\mathcal{C}=\left \{ X\in \mathbb{R}^{m\times n}: rank(X) \le r \right \}$ to convert problem (\ref{equ:1.1}) to an unconstrained optimization problem as follows:
\begin{equation}\label{equ:3.1}
\begin{aligned}
        \min_{X\in \mathbb{R}^{m\times n}} \left \| \mathcal{A}(X)-b \right \|_{2}^{2} + \delta_{\mathcal{C}}(X).
\end{aligned}
\end{equation}

Then introduce $Y\in \mathbb{R}^{m\times n}$ as an auxiliary matrix of $X$ to convert problem (\ref{equ:3.1}) as the following linear equation constrained problem:
\begin{equation}\label{equ:3.2}
    \begin{aligned}
        \min_{X,Y\in \mathbb{R}^{m\times n}} &\left \| \mathcal{A}(X)-b \right \|_{2}^{2} + \delta_{\mathcal{C}}(Y)\\
        s.t. \quad & X-Y=0.
    \end{aligned}
\end{equation}

(\ref{equ:3.2}) is the final two-block optimization problem with a linear equation constraint to describe the rank-constrained matrix sensing problem.

\section{Algorithm}\label{sec:4}

In this section, we apply ADMM to problem (\ref{equ:3.2}) to develop our algorithm for solving the rank-constrained matrix sensing problem.

The augmented Lagrangian function of problem (\ref{equ:3.2}) is
\begin{equation}\label{equ:4.1}
    \mathcal{L}_{\mu}(X,Y,\Lambda)= \left \| \mathcal{A}(X)-b \right \|_{2}^{2} + \delta_{\mathcal{C}}(Y) + \left \langle \Lambda,X-Y  \right \rangle +\frac{\mu}{2}\left \| X-Y \right \|_{F}^{2},
\end{equation}
where $\Lambda \in \mathbb{R}^{m\times n}$ is the Lagrangian multiplier of the constraint $X-Y=0$, and $\mu>0$ is the penalty parameter.

Applying ADMM to problem (\ref{equ:3.2}), the target matrix can be approximately obtained by minimizing the variables in the augmented Lagrangian function (\ref{equ:4.1}) with a Gauss-Seidel format as follows:
\begin{equation}\label{equ:4.2}
    \begin{aligned}
        \left\{\begin{matrix}
 Y^{k+1}=&\mathop{\arg\min}\limits_{Y\in \mathbb{R}^{m\times n}}\mathcal{L}_{\mu}(X^k,Y,\Lambda^{k}),\\
 X^{k+1}=&\mathop{\arg\min}\limits_{X\in \mathbb{R}^{m\times n}}\mathcal{L}_{\mu}(X,Y^{k+1},\Lambda^{k}),\\
\Lambda^{k+1}=&\Lambda^{k}+\mu (X^{k+1}-Y^{k+1}).
\end{matrix}\right.
    \end{aligned}
\end{equation}

We will solve the subproblems in Equation (\ref{equ:4.2}) one by one.

\subsection{Update Y}\label{sub:4.1}

Fix $X$ and $\Lambda$,
\begin{equation}\label{equ:4.3}
    \begin{split}
        Y^{k+1}&=\mathop{\arg\min}\limits_{Y\in \mathbb{R}^{m\times n}}\mathcal{L}_{\mu}(X^k,Y,\Lambda^{k})\\
        &=\mathop{\arg\min}\limits_{Y\in \mathbb{R}^{m\times n}}\left \{   \delta_{\mathcal{C}}(Y)+\left \langle \Lambda^{k},X^k-Y  \right \rangle +\frac{\mu}{2}\left \| X^k-Y \right \|_{F}^{2} \right \} \\
        &=\mathop{\arg\min}\limits_{Y\in \mathcal{C}}\left \| X^k-Y+\frac{\Lambda^{k}}{\mu} \right \|_{F}^{2}\\
        &=\mathcal{P}_{\mathcal{C}}(X^{k}+\frac{\Lambda^{k}}{\mu}),
    \end{split}
\end{equation}
where $\mathcal{P}_{\mathcal{C}}(\cdot)$ denotes the projection operator on the rank-constrained set $\mathcal{C}$.

In fact, the projection onto the rank-constrained set $\mathcal{C}$ has the closed solution, which is the truncated SVD \cite{TSVD,TSVD2}. First compute the SVD $X^{k}+\frac{\Lambda^{k}}{\mu}=U\Sigma V^{\top}$, then obtain $\Sigma_r$ by setting all but $r$ largest diagonal entries of $\Sigma$ to 0. And the projection is given by
\begin{equation}\label{equ:4.4}
    Y^{k+1}=\mathcal{P}_{\mathcal{C} }(X^{k}+\frac{\Lambda^{k}}{\mu})=U\Sigma_r V^{\top}.
\end{equation}

\subsection{Update X}\label{sub:4.2}

Fix $Y$ and $\Lambda$,
\begin{equation}\label{equ:4.5}
    \begin{split}
        X^{k+1}&=\mathop{\arg\min}\limits_{X\in \mathbb{R}^{m\times n}}\mathcal{L}_{\mu}(X,Y^{k+1},\Lambda^{k})\\
        &=\mathop{\arg\min}\limits_{X\in \mathbb{R}^{m\times n}}\left \{ \left \| \mathcal{A}(X)-b \right \|_{2}^{2} + \left \langle \Lambda^{k},X-Y^{k+1}\right \rangle +\frac{\mu}{2}\left \| X-Y^{k+1} \right \|_{F}^{2} \right \} \\
        &=\mathop{\arg\min}\limits_{X\in \mathbb{R}^{m\times n}}\left \{ \left \| \mathcal{A}(X)-b \right \|_{2}^{2}+\frac{\mu}{2}\left \| X-Y^{k+1}+\frac{\Lambda^{k}}{\rho} \right \|_{F}^{2} \right \}.
    \end{split}
\end{equation}

Denote $f(X)=\left \| \mathcal{A}(X)-b \right \|_{2}^{2}+\frac{\mu}{2}\left \| X-Y^{k+1}+\frac{\Lambda^{k}}{\mu} \right \|_{F}^{2}$. Obviously, it is a continuously differentiable function, hence
\begin{displaymath}
    X^{k+1}=\mathop{\arg\min}\limits_{X\in \mathbb{R}^{m\times n}}f(X) \Longleftrightarrow  \nabla f(X^{k+1})=0.
\end{displaymath}

We have 
\begin{displaymath}
    \nabla f(X^{k+1})=2\mathcal{A}^{\ast}(\mathcal{A}(X^{k+1})-b)+\mu (X^{k+1}-Y^{k+1}+\frac{\Lambda^{k}}{\mu})=(2\mathcal{A}^{\ast}\mathcal{A}+\mu I)(X^{k+1})-2\mathcal{A}^{\ast}(b)-\mu Y^{k+1}+\Lambda^{k}.
\end{displaymath}

Therefore, to update $X$, we just need to solve the matrix equation $\nabla f(X^{k+1})=0$, which is
\begin{equation}\label{equ:4.6}
    (2\mathcal{A}^{\ast}\mathcal{A}+\mu I)(X^{k+1})=2\mathcal{A}^{\ast}(b)+\mu Y^{k+1}-\Lambda^{k}.
\end{equation}

Next we solve the matrix equation (\ref{equ:4.6}).

\begin{pro}\label{pro:4.1}
The matrix equation (\ref{equ:4.6}) is equivalent to a system of linear equations, which has the unique solution given below.
\end{pro}

\begin{myproof}
Denote $\tilde{X}=vec(X^{k+1})\in \mathbb{R}^{mn}$, then $\left \langle A_i,X^{k+1} \right \rangle=vec(A_i)^{\top}\cdot vec(X^{k+1})=vec(A_i)^{\top}\cdot \tilde{X}$. And
\begin{align*}
    \mathcal{A}(X^{k+1})&=\left ( \left \langle A_1,X^{k+1} \right \rangle,\left \langle A_2,X^{k+1} \right \rangle,\cdots,\left \langle A_d,X^{k+1} \right \rangle \right )^{\top} \\
    &=\left ( vec(A_1)^{\top}\cdot \tilde{X},vec(A_2)^{\top}\cdot \tilde{X},\cdots,vec(A_d)^{\top}\cdot \tilde{X} \right )^{\top}\\ &= \left ( vec(A_1),vec(A_2),\cdots,vec(A_d) \right )^{\top}\tilde{X}.
\end{align*}
        
Denote $\tilde{A}=(vec(A_1),vec(A_2),\cdots,vec(A_d))\in \mathbb{R}^{(mn)\times d}$, then $\mathcal{A}(X^{k+1})=\tilde{A}^{\top} \tilde{X}\in \mathbb{R}^{d}$.

The adjoint operator $\mathcal{A}^{\ast}(w)=\sum_{i=1}^d w_iA_i$ for any $w\in \mathbb{R}^{d}$, and we have $vec(\mathcal{A}^{\ast}(w))=\sum_{i=1}^d w_i \cdot vec(A_i)=\left ( vec(A_1),vec(A_2),\cdots,vec(A_d) \right )w=\tilde{A}w$.

Therefore, vectorize the matrix equation (\ref{equ:4.6}). We now obtain the system of linear equations:

\begin{equation}\label{equ:4.10}
   (2\tilde{A} \tilde{A}^{\top}+\mu I_{mn})\tilde{X}=2\tilde{A}b+\mu \cdot vec(Y^{k+1})-vec(\Lambda^{k}),
\end{equation}
where $2\tilde{A}\cdot \tilde{A}^{\top}+\rho I_{mn}$ is obviously positive definite, equation (\ref{equ:4.10}) has the unique solution $\tilde{X}=(2\tilde{A} \tilde{A}^{\top}+\mu I_{mn})^{-1}(2\tilde{A}b+\mu \cdot vec(Y^{k+1})-vec(\Lambda^{k}))$.

At last, with the solution of equation (\ref{equ:4.10}) $\tilde{X}=vec(X^{k+1})$, we can reverse the vectorization operator to obtain $X^{k+1}$, which is the unique solution of the matrix equation (\ref{equ:4.6}). \hfill $\square$
\end{myproof}

\subsection{Complexity analysis}\label{sub:4.3}

Call our algorithm for solving the rank-constrained matrix sensing problem (\ref{equ:1.1}) the RCMS-ADMM, and it is summarized in Algorithm \ref{alg:1}. We will describe the parameter initialization in Section \ref{sub:7.1}.

\begin{algorithm}[H]
  \KwIn{The measurement vector $b$, the linear operator $\mathcal{A}$ and an upper estimation of rank $r$;}
  \KwOut{The approximate low-rank matrix $\Hat{X}$;}
  initialization:$X^{0},\Lambda^{0},\mu$, and let $k=0$;
  
    \Do{not satisfy the stop criterion}{
      Update $Y$ by $Y^{k+1}=\mathcal{P}_{\mathcal{C} }(X^{k}+\frac{\Lambda^{k}}{\mu})$;
      
      Update $X$ by solving the matrix equation (\ref{equ:4.6});
      
      Update $\Lambda$ by $\Lambda^{k+1}=\Lambda^{k}+\mu (X^{k+1}-Y^{k+1})$;
      
      Update $k$ by $k=k+1$;
    }
    
    \Return $\Hat{X}=Y^{k}$.
  \caption{RCMS-ADMM}
  \label{alg:1}
\end{algorithm}

For the general linear operator $\mathcal{A}$, $A_i (1\le i\le d)$ are linearly independent without any special structure. Therefore, when updating $X$, we need to solve the system of linear equations (\ref{equ:4.10}) with the order of coefficient matrix $mn$. Note that
$2\tilde{A} \tilde{A}^{\top}+\mu I_{mn}$ does not change in the iteration; therefore, the inverse $(2\tilde{A} \tilde{A}^{\top}+\mu I_{mn})^{-1}$ can be calculated only once outside the iteration. But it will still cost a lot of time with $mn$ increasing. An obvious way to accelerate the calculation is using the SMW-formula \cite{SMW1,SMW2}, converting the inverse of an $(mn)\times (mn)$ matrix to the inverse of a $d\times d$ matrix:
\begin{equation}\label{equ:4.13}
    (2\tilde{A} \tilde{A}^{\top}+\mu I_{mn})^{-1}=\frac{I_{mn}}{\mu}-\frac{2\tilde{A}(I_d+\frac{2\tilde{A}^{\top} \tilde{A}}{\mu})^{-1}\tilde{A}^{\top}}{\mu^{2}}.
\end{equation}

In each iteration, updating $Y$ by truncated SVD of an $m\times n$ matrix (\ref{equ:4.4}) costs $\mathcal{O}(r^2\cdot \min(m,n))$ time. Since $(2\tilde{A}\cdot \tilde{A}^{\top}+\rho I_{mn})^{-1}$ is calculated only once outside the iteration, updating $X$ by solving (\ref{equ:4.10}) only involves the multiplication between an $(mn)\times (mn)$ matrix and an $mn$ vector, which will cost $\mathcal{O}(m^2 n^2)$ time. Therefore, when the linear operator $\mathcal{A}$ has no special structure, the time complexity of Algorithm \ref{alg:1} is $\mathcal{O}(m^2 n^2)$ each iteration.

\section{A specific case: rank-constrained matrix completion problem}\label{sec:5}

In this section, we consider a special but significant application of rank-constrained matrix sensing, which is known as the rank-constrained matrix completion problem. In this specific case, each matrix $A_k=(a_{ij}^k)_{m\times n}, k=1,2,\cdots,d$ in the linear operator 
\begin{align*}
 \mathcal{A}:&\mathbb{R}^{m\times n}\to \mathbb{R}^d  \\
&X\mapsto (\left \langle A_1,X \right \rangle,\left \langle A_2,X \right \rangle,\cdots,\left \langle A_d,X \right \rangle)^{\top}
\end{align*}
has exactly one non-zero entry, and $\mathcal{A}$ actually returns a subset of the target matrix.

Denote the location matrix $\Omega \in \mathbb{R}^{m\times n}$, where
\begin{displaymath}
    \Omega_{ij}=\begin{cases}
     1, & \text{ \textit{if} } a_{ij}^k\ne 0 \text{ \textit{for some} } 1\le k\le d \\
     0, & \text{ \textit{otherwise} } 
    \end{cases}
\end{displaymath}
Let $M=\Omega \odot X$. There is a one-to-one correspondence between $b=\mathcal{A}(X)$ and $M_{ij}$ where $(i,j)$ satisfies $\Omega_{ij}=1$, which means $M$ is the measurement matrix in this case. Therefore, $\mathcal{A}^{\ast}(b)=\sum_{i=1}^d b_i A_i=M$, and $2\mathcal{A}^{\ast}\mathcal{A}(X^{k+1})=2\sum_{i=1}^d \left \langle A_i,X^{k+1} \right \rangle A_i=2\Omega\odot X^{k+1}$.  

The matrix equation (\ref{equ:4.6}) becomes
\begin{displaymath}
    (2\Omega+\mu \mathbf{1})\odot X^{k+1}=2M+\mu Y^{k+1}-\Lambda^k,
\end{displaymath}
whose solution can be directly given by hadamard division:
\begin{equation}\label{equ:5.1}
    X^{k+1}=(2M+\mu Y^{k+1}-\Lambda^k)\oslash (2\Omega+\mu \mathbf{1}).
\end{equation}

The algorithm for the specific rank-constrained matrix completion problem is summarized in the following.

\begin{algorithm}[H]
  \KwIn{The measurement matrix $M$, the location matrix $\Omega$, and an upper estimation of rank $r$;}
  \KwOut{The approximate low-rank matrix $\Hat{X}$;}
  initialization:$X^{0},\Lambda^{0},\mu$, and let $k=0$;
  
    \Do{not satisfy the stop criterion}{
      Update $Y$ by $Y^{k+1}=\mathcal{P}_{\mathcal{C} }(X^{k}+\frac{\Lambda^{k}}{\mu})$;
      
      Update $X$ by $X^{k+1}=(2M+\mu Y^{k+1}-\Lambda^k)\oslash (2\Omega+\mu \mathbf{1})$;
      
      Update $\Lambda$ by $\Lambda^{k+1}=\Lambda^{k}+\mu (X^{k+1}-Y^{k+1})$;
      
      Update $k$ by $k=k+1$;
    }
    
    \Return $\Hat{X}=Y^{k}$.
  \caption{RCMC-ADMM}
  \label{alg:2}
\end{algorithm}

Unlike the Algorithm \ref{alg:1}, updating $X$ in Algorithm \ref{alg:2} only involves a hadamard division, which costs $\mathcal{O}(mn)$ time. Therefore, the time complexity of Algorithm \ref{alg:2} is $\mathcal{O}(mn+r^2\cdot \min(m,n))$ per iteration.

\section{Convergence analysis}\label{sec:6}

In this section, we will establish the global convergence of the nonconvex Algorithm \ref{alg:1} based on the KL property of functions in problem (\ref{equ:3.2}).

\begin{pro}\label{pro:6.1}
    $f(X)=\left \| \mathcal{A}(X)-b \right \|_2^2$ is a Gradient-$L$-Lipschitz continuous and KL function. 
\end{pro}

\begin{myproof}
    Since $f(X)=\left \| \mathcal{A}(X)-b \right \|_2^2=\sum_{i=1}^d (\left \langle A_i,X \right \rangle -b_i)^2$ is a real polynomial, it is of course a KL function.

    $\nabla f(X)=2\mathcal{A}^{\ast}(\mathcal{A}(X)-b)=2\sum_{i=1}^d (\left \langle A_i,X \right \rangle -b_i)\cdot A_i$, therefore we have
    \begin{equation}\label{equ:6.1}
    \begin{split}
        \left \| \nabla f(X)-\nabla f(Y) \right \|_F &= \left \| 2\sum_{i=1}^d \left \langle A_i,X-Y \right \rangle \cdot A_i \right \|_F \\
        &\le 2\sum_{i=1}^d \left \| \left \langle A_i,X-Y \right \rangle \cdot A_i \right \|_F \\
        &= 2\sum_{i=1}^d \left | \left \langle A_i,X-Y \right \rangle \right | \cdot \left \|   A_i \right \|_F \\
        &\le (2\sum_{i=1}^d \left \|   A_i \right \|_F^2)\cdot \left \|   X-Y \right \|_F,
    \end{split}
    \end{equation}
    where the first inequality is due to the triangular inequality; the second equality is because of the norm being positively homogeneous; and the last inequality uses the Cauchy-Schwarz inequality. 

    From the inequality (\ref{equ:6.1}), we know that $f(X)$ is a Gradient-$L$-Lipschitz continuous function with the Lipschitz constant $L=2\sum_{i=1}^d \left \|   A_i \right \|_F^2>0$.
    \hfill $\square $
\end{myproof}

\begin{pro}[\cite{rKL}]\label{pro:6.2}
    The rank function $f(X)=rank(X)$ is proper and lower semicontinuous. And it is also a semi-algebraic, therefore a KL function. 
\end{pro}

Due to Proposition \ref{pro:2.1}, the next Corollary \ref{cor:6.1} holds immediately.

\begin{cor}\label{cor:6.1}
    $f(X)=rank(X)$ is a closed function.
\end{cor}

\begin{cor}\label{cor:6.2}
    The indicator function 
    \begin{displaymath}
    \delta_{\mathcal{C}}(X)=\begin{cases}
    0, & \text{ \textit{if} } X\in \mathcal{C} \\
    +\infty, & \text{ \textit{if} } X\notin \mathcal{C}
    \end{cases}
    \end{displaymath}
    to the rank-constrained set
    $\mathcal{C}=\left \{ X\in \mathbb{R}^{m\times n}:rank(X)\le r \right \}$ is proper and lower semicontinuous, and it is also a KL function.
\end{cor}

\begin{myproof}
    From the definition of the proper function (\ref{def:2.3}), it is easy to see that $\delta_{\mathcal{C}}(X)$ is proper.

    Next we prove that $\delta_{\mathcal{C}}(X)$ is closed, therefore it is lower semicontinuous.

    According to the definition of the closed function (\ref{def:2.3}), we only need to prove that for any sequence $(X_{k},t_{k})\in epi\delta_{\mathcal{C}}$, if $(X_{k},t_{k})\to (X,t)$, then $t\ge \delta_{\mathcal{C}}(X)$.

    Since $t_{k}\to t$, there are at most the following two cases when $k$ is large enough:

    If $rank(X_k)>r$, then $t_{k}\ge \delta_{\mathcal{C}}(X_{k})=+\infty$. Therefore, $t=\lim_{k \to +\infty} t_k = +\infty \ge \delta_{\mathcal{C}}(X)$.

    If $rank(X_k)\le r$, then $t_{k}\ge \delta_{\mathcal{C}}(X_{k})=0$. Hence $t=\lim_{k \to +\infty} t_k \ge 0$. Corollary \ref{cor:6.1} tells that the rank function is closed. Therefore, 
    \begin{center}
        $\left.
        \begin{aligned}
        &(X_{k},r)\in epirank(\cdot)  \\
        &(X_{k},r)\to (X,r)
        \end{aligned}
        \right\}
        \Longrightarrow
        (X,r)\in epirank(\cdot)$,
    \end{center}
    which means $r\ge rank(X)$. Hence
    $t\ge 0= \delta_{\mathcal{C}}(X)$.

    Above all, $t\ge \delta_{\mathcal{C}}(X)$, which tells that $\delta_{\mathcal{C}}(X)$ is closed.

    Then we illustrate that $\mathcal{C}$ is semi-algebraic. 
    
    Since $rank(\cdot)\in \mathbb{N}$, we can express $\mathcal{C}$ as
    \begin{displaymath}
        \mathcal{M}=\left \{ X\in \mathbb{R}^{m\times n}:rank(X)\le r \right \}=\bigcup_{i=0}^{r}\left \{ X\in \mathbb{R}^{m\times n}:rank(X)=i \right \}.
    \end{displaymath}

    Proposition \ref{pro:6.2} tells that $rank(\cdot)$ is semi-algebraic. Therefore, through the definition of the semi-algebra (\ref{def:2.4}), $graphrank(\cdot)=\left \{ (X,i):rank(X)=i \right \}$
    is a semi-algebraic set.

    Denote the projection operator
    \begin{align*}
        \pi:&\mathbb{R}^{m\times n}\times \mathbb{R}\to \mathbb{R}^{m\times n}  \\
        &(X,i)\mapsto X,
    \end{align*}
    Then $\mathcal{C}=\bigcup_{i=0}^{r}\pi(graphrank(\cdot))$ is also semi-algebraic due to its stability under the finite union and projection. 

    At last, since $\delta_{\mathcal{C}}(X)$ is the indicator function to the semi-algebra $\mathcal{C}$, it is also a semi-algebra. From the Proposition \ref{pro:2.3}, $\delta_{\mathcal{C}}(X)$ is a KL function.
    \hfill $\square$
\end{myproof}

In the following part of the section, we refer to \cite{converge} to establish the global convergence of the Algorithm \ref{alg:1}.

Consider the following class of the two-block nonconvex optimization problem with a linear equality constraint
\begin{equation}\label{equ:6.2}
    \begin{aligned}
        \min_{x,y} \enspace &f(x)+g(y)\\
        s.t.\enspace &Ax+y=b
    \end{aligned}
\end{equation}
with the following four assumptions.
\begin{ass}[\cite{converge}]\label{ass:6.1}
    There are four assumptions for problem (\ref{equ:6.2}):
    \begin{enumerate}[(i)]
        \item $f$ is Gradient-$L$-Lipschitz continuous and $g$ is proper lower semicontinuous,
        \item The subproblems in (\ref{equ:4.2}) has solutions,
        \item The penalty parameter $\mu>2L$,
        \item  $A^{\top}A\succeq \rho I_n$ for some $\rho>0$.
\end{enumerate}
\end{ass}

Let $f(X)=\left \| \mathcal{A}(X)-b  \right \|_2^2$, $g(Y)=\delta_{\mathcal{C}}(Y)$, $A=-I_m$, $b=0$ in problem (\ref{equ:6.2}), we obtain the problem (\ref{equ:3.2}). Based on the Proposition \ref{pro:6.1}, Corollary \ref{cor:6.2} and the analysis in Section \ref{sec:4}, it is apparent that problem (\ref{equ:3.2}) satisfies the Assumption \ref{ass:6.1}.

The main convergence result about problem (\ref{equ:6.2}) in \cite{converge} is shown below.
\begin{lem}[\cite{converge}]\label{lem:6.1}
    Applying the classic ADMM to problem (\ref{equ:6.2}), we obtain a sequence 
    $\left \{w^{k}=(x^k,y^k,\lambda^k)\right \}$. If 
    \begin{equation}\label{equ:6.3}
    \bar{f}:=\inf_{x}\left \{ f(x)-\frac{1}{2L}\left \| \nabla f(x) \right \|^{2} \right \}>-\infty,
    \end{equation}
    $\liminf_{\left \|x  \right \|\to +\infty}f(x)=+\infty$ and 
    $\inf_{y}g(y)>-\infty$, then 
    $\left \{w^{k}\right \}$ is bounded. 
    Besides, if $f$ and $g$ are also KL functions, then $\left \{ w^{k} \right \}$ converges to the KKT point of (\ref{equ:6.2}).  
\end{lem}

Based on the above Lemma (\ref{lem:6.1}), we can conclude that the Algorithm \ref{alg:1} is globally convergent.

\begin{thm}\label{thm:6.1}
    Algorithm \ref{alg:1} is globally convergent, and the sequence $\left \{ (Y^k,X^k,\Lambda^k) \right \}$ generated by it globally converges to the KKT point of the problem (\ref{equ:3.2}).
\end{thm}

\begin{myproof}
    For $f(X)=\left \| \mathcal{A}(X)-b  \right \|_2^2$, we know $\nabla f(X)=2\mathcal{A}^{\ast}(\mathcal{A}(X)-b)$ and $L=2\sum_{i=1}^d \left \| A_i  \right \|_F^2$ from Proposition \ref{pro:6.1}.
    Therefore, for any $X\in \mathbb{R}^{m\times n}$,
    \begin{equation}\label{equ:6.4}
    \begin{split}
        \left \| \mathcal{A}(X)-b  \right \|_2^2-\frac{1}{2L}\left \| 2\mathcal{A}^{\ast}(\mathcal{A}(X)-b)  \right \|_F^2 &= \sum_{i=1}^d (\left \langle A_i,X \right \rangle-b_i)^2 - \frac{1}{2L}\left \| 2\sum_{i=1}^d (\left \langle A_i,X \right \rangle-b_i)\cdot A_i  \right \|_F^2  \\
        &\ge \sum_{i=1}^d (\left \langle A_i,X \right \rangle-b_i)^2-\frac{1}{2L}\cdot 4 \sum_{i=1}^d (\left \langle A_i,X \right \rangle-b_i)^2\cdot \left \| A_i  \right \|_F^2\\    
        &\ge \sum_{i=1}^d (\left \langle A_i,X \right \rangle-b_i)^2-\frac{2}{L}(\sum_{i=1}^d (\left \langle A_i,X \right \rangle-b_i)^2)(\sum_{i=1}^d \left \| A_i  \right \|_F^2) \\
        &= 0,
    \end{split}
    \end{equation}
    where the first inequality is due to the triangular inequality; and the second inequality is because $\sum_{i=1}^n a_i^2 b_i^2 \le (\sum_{i=1}^n a_i^2) (\sum_{i=1}^n b_i^2)$ for any $a_i,b_i\in \mathbb{R}$.

    Hence (\ref{equ:6.3}) in Lemma \ref{lem:6.1} is satisfied. 
    
    It is easy to see that $\liminf_{\left \|X  \right \|\to +\infty}f(X)=\lim_{\left \|X  \right \|\to +\infty}\sum_{i=1}^d (\left \langle A_i,X \right \rangle-b_i)^2 =+\infty$ and 
    $\inf_{Y}g(Y)=\inf_{Y}\delta_{\mathcal{C}}(Y)=0>-\infty$. Besides, according to Proposition \ref{pro:6.1} and Corollary \ref{cor:6.2} we know $f$ and $g$ are both KL functions. Therefore, by Lemma \ref{lem:6.1}, the generated sequence $\left \{ (Y^k,X^k,\Lambda^k) \right \}$ is globally convergent to the KKT point of the problem (\ref{equ:3.2}).
    \hfill $\square$
\end{myproof}

Moreover, if the original matrix $X$ is low-rank such that the rank-constrained matrix sensing problem (\ref{equ:1.1}) has the unique solution \footnote{For the general linear operator $\mathcal{A}$, we only know that problem (\ref{equ:1.1}) has the unique solution when the rank of $X$ is low, but the precise relation between the uniqueness and rank of $X$ is still an open problem. However, for the specific rank-constrained matrix completion problem, $d=\mathcal{O}(n^{1.2}r\log_{10}n)$ known entries are enough to recover the original $n\times n$ rank $r$ matrix \cite{lrcc}.}, then we have the following theorem to verify if Algorithm \ref{alg:1} converges to it during the iterations.

\begin{thm}\label{thm:6.2}
     Suppose that the rank-constrained matrix sensing problem (\ref{equ:1.1}) has the unique solution, then the KKT point which $\left \{ (Y^k,X^k,\Lambda^k) \right \}$ generated by Algorithm \ref{alg:1} converges to is actually the unique solution of (\ref{equ:3.2}) \textit{iff} the Lagrangian multiplier $\Lambda^k \to 0$. 
\end{thm}

\begin{myproof}
    If $(Y^k,X^k,\Lambda^k)\to (\bar{Y},\bar{X},\bar{\Lambda})$ which is the unique solution of problem (\ref{equ:3.2}), then it must also be the KKT point of problem (\ref{equ:3.2}) due to the linear equality constraint \cite{cop}. The Lagrangian function of problem (\ref{equ:3.2}) is 
    \begin{equation}\label{equ:6.5}
        \mathcal{L}(Y,X,\Lambda)=\left \| \mathcal{A}(X)-b  \right \|_2^2+\delta_{\mathcal{C}}(Y)+\left \langle \Lambda,X-Y \right \rangle,
    \end{equation}
    and it satisfies the KKT conditions at $(\bar{Y},\bar{X},\bar{\Lambda})$, which are
    \begin{subnumcases}
    {}
        0 = \nabla_{X} \mathcal{L}(\bar{Y},\bar{X},\bar{\Lambda})=2\mathcal{A}^{\ast}(\mathcal{A}(\bar{X})-b)+\bar{\Lambda}, \label{equ:6.6a} \\
        0 \in \partial_{Y}\mathcal{L}(\bar{Y},\bar{X},\bar{L})=\partial_{Y}\delta_{\mathcal{C}}(\bar{Y})-\bar{\Lambda}, \label{equ:6.6b} \\
        0 =\bar{X}-\bar{Y}. \label{equ:6.6c}
    \end{subnumcases}

    From the updating rule of $Y$ in (\ref{equ:4.3}), we know that $\delta_{\mathcal{C}}(Y^k)=0$ always holds in the iterations. Since $\delta_{\mathcal{C}}(Y)$ is proper lower semicontinuous by Proposition \ref{pro:6.2}, we obtain $0\le \delta_{\mathcal{C}}(\bar{Y})\le \liminf_{k\to +\infty}\delta_{\mathcal{C}}(Y^{k})=0$. Hence $\delta_{\mathcal{C}}(Y)\equiv 0$ when limiting the space on the generated points and their accumulations by Algorithm \ref{alg:1}. Besides, the upper estimation rank $r$ in problem (\ref{equ:3.2}) is always set larger than the real rank of the unique solution $\bar{X}$ in practice \footnote{The solution of problem (\ref{equ:1.1}) is unique as long as $\text{rank}(\bar{X})\le r\le \bar{r}$, where $\bar{r}$ is the precise upper bound of rank(X) for the uniqueness of solutions in theory (though we cannot know it for the general linear operator $\mathcal{A}$). 
    In practice, the rank of the unique solution $\bar{X}$ is always low. Hence we prefer to set the upper estimation of the rank $r$ larger than the real rank($\bar{X}$).}. Therefore, $\partial_{Y}\delta_{\mathcal{C}}(\bar{Y})=\widetilde{\mathcal{N}}_{\mathcal{C}}(\bar{Y})=\{0\}$. From (\ref{equ:6.6b}) we obtain $\bar{\Lambda}=0$, which proves $\Lambda^k \to 0$.

    As for the sufficiency, Since the limit point $(\bar{Y},\bar{X},\bar{\Lambda})$ is the KKT point of problem (\ref{equ:3.2}), it must satisfy (\ref{equ:6.6a})-(\ref{equ:6.6c}). Thus we have $\mathcal{A}^{\ast}(\mathcal{A}(\bar{X})-b)=0$ from (\ref{equ:6.6a}) since $\Lambda^k\to0$. Because $\mathcal{A}^{\ast}(w)=\sum_{i=1}^d w_i A_i$ for $w\in \mathbb{R}^d$, if $\mathcal{A}^{\ast}(w)=0$, there must be $w=0$. Otherwise, there exists some $A_k$ which can be expressed by other $A_i,i\ne k$. Then the measurement vector $b_k=\left \langle A_k,X \right \rangle$ can also be expressed by other measurements $b_i,i\ne k$, which suggests that this measurement is redundant. Hence from $\mathcal{A}^{\ast}(\mathcal{A}(\bar{X})-b)=0$ we have $\mathcal{A}(\bar{X})=b$. From (\ref{equ:6.6c}) we know $\bar{X}=\bar{Y}$, which means that $rank(\bar{X})=rank(\bar{Y})\le r$ since 
    we have proved that $\delta_{\mathcal{C}}(\bar{Y})=0$. Therefore, $\bar{X}$ also satisfies $\mathcal{A}(\bar{X})=b$, $rank(\bar{X})\le r$, which means that $\bar{X}$ is actually the unique optimal solution of problem (\ref{equ:1.1}).
    \hfill $\square$
\end{myproof}

\begin{rem}
    Based on the Theorem \ref{thm:6.2}, we can check if $\left\| \Lambda^k \right\|_{F}<Tol$ as an extra stopping criterion. Therefore, when Algorithm \ref{alg:1} stops earlier, we know that the output $X^k$ is actually the unique optimal solution of problem (\ref{equ:1.1}). Besides, for rank-constrained matrix sensing problem, it has no spurious local minima when the measurements $\{A_i\}$ satisfy mild conditions \cite{nolocalmin}; and \cite{strictsaddle} points out that the strict saddle points set of lots of first-order methods has measure zero. Hence the Algorithm \ref{alg:1} can stop earlier to find the unique optimal solution most of the time in practice.
\end{rem}

\section{Numerical simulations}\label{sec:7}

In this section, we consider the specific rank-constrained matrix completion problem, and design numerical experiments to verify the proposed specialized Algorithm \ref{alg:2} by comparing it with the other algorithms. 

\subsection{Experiment setup}\label{sub:7.1}

The experiments are performed on a PC with an Intel Core i7-13700H 2.4 GHz CPU and 16 GB RAM. The results in Sections \ref{sub:7.4}-\ref{sub:7.6} have been averaged over 10 trials. The original matrix $X\in \mathbb{R}^{m\times n}$ is generated by $X=BC^\top$ where $B\in \mathbb{R}^{m\times r}$ and $C\in \mathbb{R}^{n\times r}$ both have the $i.i.d.$ Gaussian random entries from $\mathcal{N}(0,1)$. And the measurement $b$ is $d$ randomly chosen from $X$, which might be added with an extra Gaussian white noise. The recovered matrix is denoted as $\hat{X}$. And we use $\text{SNR}_{\text{r}}=20\log_{10}(\|X\|_F/\|X-\widehat{X}\|_F)$ and $\text{SNR}_{\text{m}}=20\log_{10}(\|b\|_2/\|e\|_2)$ to measure the reconstruction error and measurement noise level, respectively. The computational cost is measured by the number of iterations and running time.

Because the projection in (\ref{equ:4.3}) is a hard thresholding operator, we check the NIHT \cite{NIHT} for the RCMC. The NN-ADMM \cite{NNADMM} is the same ADMM scheme algorithm for the nuclear norm regularized least squares problem. We introduce it in our experiments as the most direct comparison and call it NN-ADMM for distinguishing with our Algorithm \ref{alg:2} called RC-ADMM for simplicity.

The parameter settings are shown in Table \ref{tab:1}. For Algorithm \ref{alg:2}, $X^{0}$ and $L^{0}$ are initialized as an $i.i.d.$ Gaussian random matrix and a null matrix, respectively.
The maximum iterations MaxIter are all 500; and the basic stop criterion is $\|X^{k+1}-X^k\|_F/\|X^k\|_F<\text{Tol}$ with Tol$=10^{-4}$.   

\begin{table}[hbt!]
    \centering
    \begin{tabular}{cc}
    \hline
        Algorithm & Parameters \\\hline
        NIHT \cite{NIHT} & $r$, Tol$=10^{-4}$, MaxIter=500 \\
        NN-ADMM \cite{NNADMM} & step size $\mu=10^{-4}$, ratio $\rho=1.1$, Tol$=10^{-4}$, MaxIter=500 \\
        RC-ADMM & $r$, $\mu=1$, Tol$=10^{-4}$, MaxIter=500 \\ \hline 
    \end{tabular}
    \caption{Parameter settings}
    \label{tab:1}
\end{table}

\subsection{Validity to recover from incomplete observations}\label{sub:7.2}

\begin{figure*}[!t]
	\centering
	\includegraphics[width=5in]{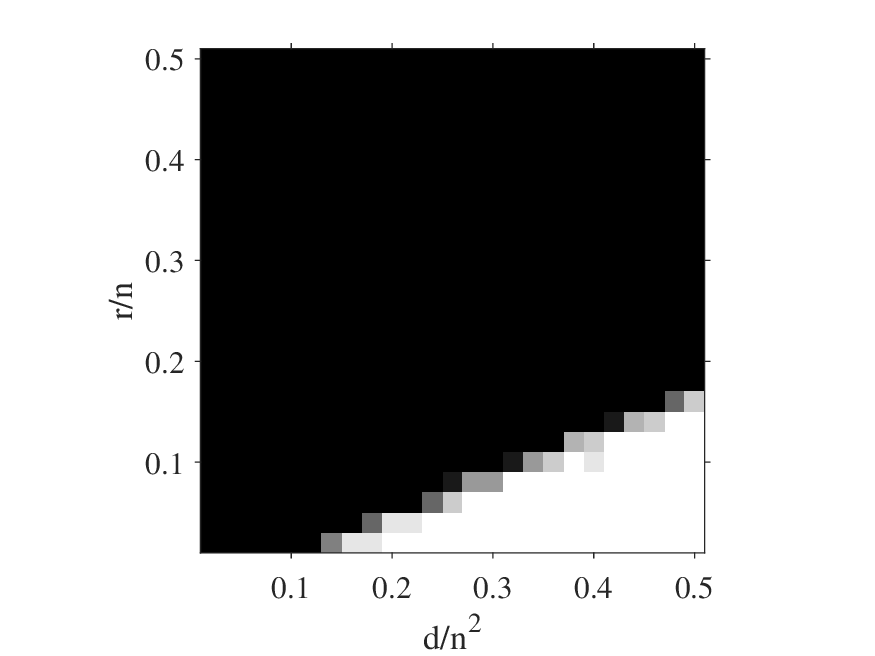}
	\caption{Phase transition of matrix completion problem, $m=n=100$.}
	\label{fig:1}
\end{figure*}

Firstly, we verify that the Algorithm \ref{alg:2} has the ability to recover the original matrix from incomplete observations with varying the rank from varying the sampling rate. Let $m=n=100$, the sampling rate $d/n^2=[0.02:0.02:0.5]$, the rank $r/n=[0.02:0.02:0.5]$. We take 10 trials for each $(\frac{r}{n},\frac{d}{n^2})$-pair without stop criterion, and say a trial to be successful if $\text{SNR}_{\text{r}}\ge 70$ dB (black = 0\% and white = 100\%). The phase transition result is shown in Figure \ref{fig:1}. It can be seen that when the rank of the original matrix is low, Algorithm \ref{alg:2} can reconstruct it with a high precision, and the sampling rate can be really low at this time.

\subsection{Test for the Lagrangian multiplier}\label{sub:7.3}

Theorem \ref{thm:6.2} tells that the Algorithm \ref{alg:1} and \ref{alg:2} output an optimal solution \textit{iff} the Lagrangian multiplier $\Lambda^k\to 0$ when considering the noiseless case. We check it in this experiment. Set $m=n=500$, $r=20$, $d/mn=0.2$, and observe the average change of $\Lambda^k$ in iterations over 10 random trials.  

\begin{figure*}[!t]
	\centering
	\subfloat[]{\includegraphics[width=2.0in]{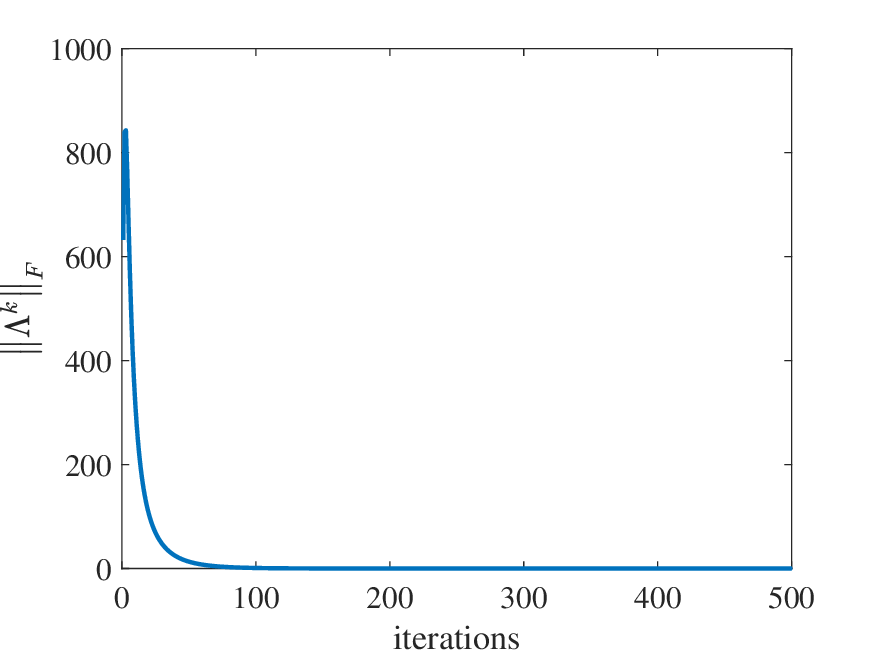}}
	\hfil
	\subfloat[]{\includegraphics[width=2.0in]{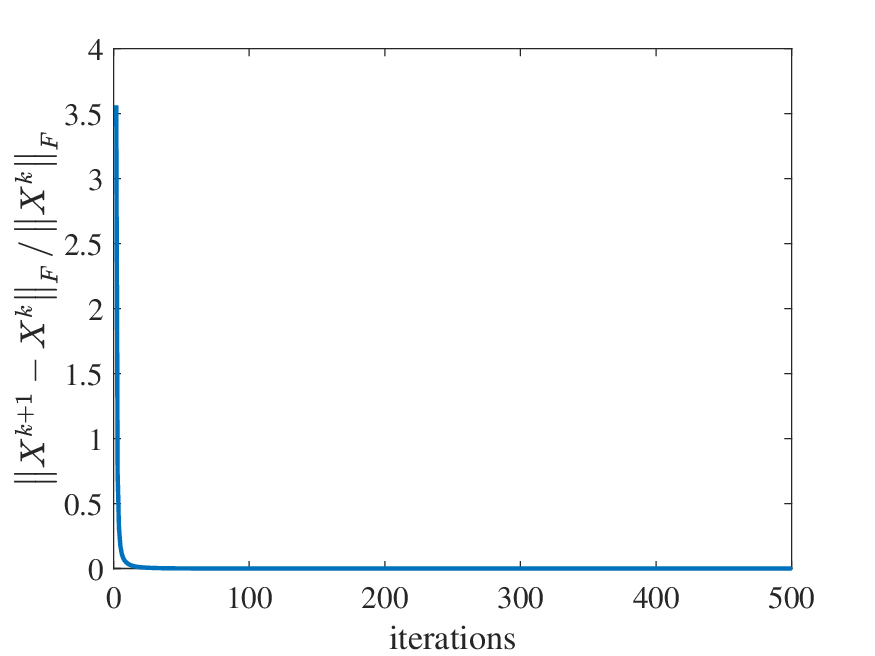}}
        \hfil
        \subfloat[]{\includegraphics[width=2.0in]{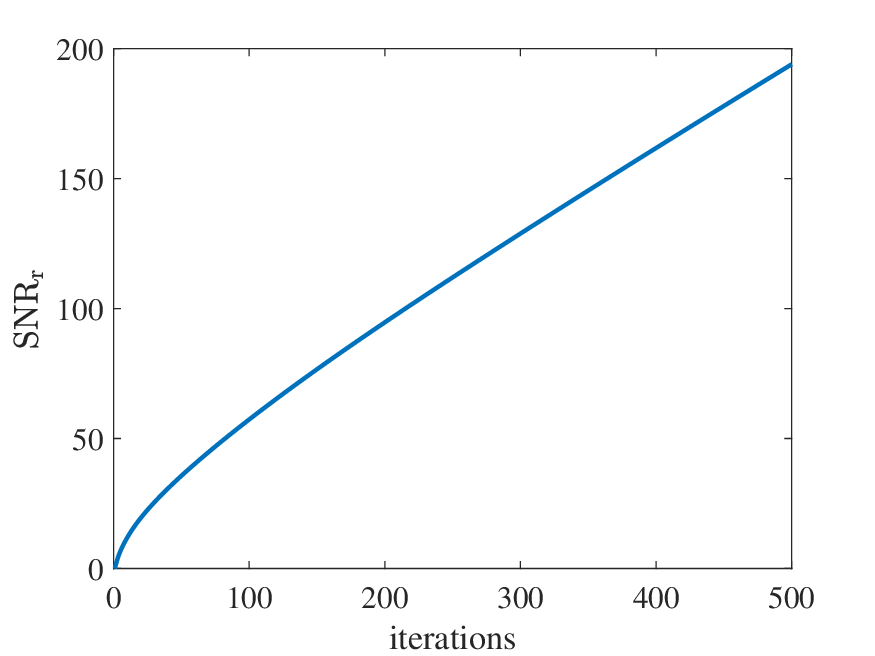}}
	\caption{The average change in iterations over 10 experiments, $m=n=500$, $r=20$, $d/mn=0.2$. (a) The average change of the Lagrangian multiplier $\Lambda^k$. (b) The average change of the relative error $\|X^{k+1}-X^k\|_F/\|X^k\|_F$. (c) The average change of the recovery precision $\text{SNR}_{\text{r}}$.}
	\label{fig:2}
\end{figure*}

As shown in Figure \ref{fig:2}, the Lagrangian multiplier indeed converges to zero in iterations, which is consistent with the theorem \ref{thm:6.2}. Hence when considering noiseless case, we can use the Lagrangian multiplier to see if the limit point is an optimal solution.

\subsection{Performance with different sampling rate}\label{sub:7.4}

Starting with this experiment, we introduce the same stop criterion $\|X^{k+1}-X^k\|_F/\|X^k\|_F<\text{Tol}$ with Tol$=10^{-4}$ for all the tested algorithms. Since there is always noise when observing in reality, we also assume that the measurement noise is $\text{SNR}_{\text{m}}=20$ in the following three experiments. 

In this section, we study how does the sampling rate $d/mn$ affect the reconstruction. Set $m=n=500$, $r=10$, and let the sampling rate $d/mn=[0.06:0.02:0.24]$. The results are shown in Table \ref{tab:2}.

\begin{table}[hbt!]
    \centering
        \resizebox{\textwidth}{!}{
    \begin{tabular}{c|ccc|ccc|ccc}
 \hline
\multirow{2}*{$d/mn$} & \multicolumn{3}{c|}{$\text{SNR}_{\text{r}}$} & \multicolumn{3}{c|}{Iter} & \multicolumn{3}{c}{Time (seconds)} \\
           & NIHT & NN-ADMM & RC-ADMM & NIHT & NN-ADMM & RC-ADMM & NIHT & NN-ADMM & RC-ADMM \\ \hline
           0.06 & 6.56 & 4.70 & 13.45 & 477 & 95 & 500 & 1.63 & 3.19 & 2.24 \\
           0.08 & 15.94 & 8.47 & 19.33 & 231 & 94 & 298 & 0.84 & 3.13 & 1.57 \\
           0.10 & 19.65 & 12.06 & 21.30 & 127 & 93 & 178 & 0.51 & 3.14 & 0.97 \\
           0.12 & 18.31 & 14.79 & 22.56 & 151 & 94 & 122 & 0.56 & 3.17 & 0.68 \\
           0.14 & 21.23 & 16.58 & 23.61 & 91 & 92 & 91 & 0.37 & 3.11 & 0.53 \\
           0.16 & 24.37 & 17.68 & 24.37 & 38 & 89 & 71 & 0.19 & 3.01 & 0.40 \\
           0.18 & 25.04 & 18.41 & 25.04 & 31 & 87 & 57 & 0.16 & 2.92 & 0.33 \\
           0.20 & 25.58 & 19.13 & 25.58 & 29 & 86 & 48 & 0.14 & 2.90 & 0.27 \\
           0.22 & 26.05 & 19.58 & 26.05 & 26 & 85 & 40 & 0.13 & 2.89 & 0.23 \\
           0.24 & 26.48 & 20.01 & 26.48 & 24 & 85 & 34 & 0.12 & 2.88 & 0.20 \\
    \hline
    \end{tabular}
    }
    \caption{Comparison of the algorithms with different $d$: $\text{SNR}_{\text{m}}=20$, $m=n=500$, $r=10$}
    \label{tab:2}
\end{table}

It can be seen that both NIHT and RC-ADMM perform good when the sampling rate is not so low ($i.e.$, $d/mn \ge 0.16$), while there is a clear gap between NN-ADMM and them in both recovery accuracy and speed. NIHT is the most efficient, and RC-ADMM is also fast enough to reconstruct the matrix. However, when the sampling rate is really low ($i.e.$, $d/mn \le 0.14$), RC-ADMM demonstrates a clear advantage in recovery precision than NIHT, which is the key of matrix completion.

\subsection{Performance with different rank}\label{sub:7.5}

The fifth experiment studies how does the rank $r$ influence the recovery. Set $m=n=500$, $d/mn=0.20$ and test the rank $r=[2:4:38]$. 

\begin{table}[hbt!]
    \centering
   \resizebox{\textwidth}{!}{
    \begin{tabular}{c|ccc|ccc|ccc}
 \hline
\multirow{2}*{$r$} & \multicolumn{3}{c|}{$\text{SNR}_{\text{r}}$} & \multicolumn{3}{c|}{Iter} & \multicolumn{3}{c}{Time (seconds)} \\
           & NIHT & NN-ADMM & RC-ADMM & NIHT & NN-ADMM & RC-ADMM & NIHT & NN-ADMM & RC-ADMM \\ \hline
           2 & 33.36 & 23.33 & 33.36 & 13 & 89 & 30 & 0.05 & 2.98 & 0.14 \\
           6 & 28.18 & 21.12 & 28.18 & 21 & 87 & 35 & 0.10 & 2.95 & 0.20 \\
           10 & 25.64 & 19.12 & 25.64 & 28 & 86 & 48 & 0.14 & 2.91 & 0.27 \\
           14 & 23.67 & 17.05 & 23.67 & 39 & 86 & 63 & 0.22 & 2.90 & 0.40 \\
           18 & 22.07 & 14.78 & 22.07 & 53 & 87 & 83 & 0.33 & 2.91 & 0.57 \\
           22 & 20.61 & 12.21 & 20.62 & 73 & 85 & 110 & 0.55 & 2.84 & 0.87 \\
           26 & 19.14 & 9.73 & 19.15 & 108 & 82 & 148 & 0.82 & 2.74 & 1.18 \\
           30 & 17.56 & 7.64 & 17.59 & 157 & 80 & 214 & 1.24 & 2.68 & 1.74 \\
           34 & 15.79 & 6.24 & 15.84 & 264 & 79 & 338 & 2.27 & 2.62 & 3.20 \\
           38 & 13.41 & 5.16 & 13.63 & 500 & 78 & 500 & 4.68 & 2.56 & 5.21 \\
    \hline
    \end{tabular}
    }
    \caption{Comparison of the algorithms with different $r$: $\text{SNR}_{\text{m}}=20$, $m=n=500$, $d/mn=0.2$}
    \label{tab:3}
\end{table}

As shown in Table \ref{tab:3}, The recovery accuracy of NIHT and RC-ADMM is almost the same for the tested $r$. And the recovery accuracy of RC-ADMM is only slightly larger than that of NIHT when the rank $r$ is large ($i.e.$, $r\ge 30$). It is apparent that NN-ADMM performs worse than NIHT and RC-ADMM in recovery precision. 
This experiment implies that RC-ADMM performs better than the other two algorithms even if the original matrix is not really low-rank.

\subsection{Performance with different matrix size}\label{sub:7.6}

In the last experiment, we test the performance of the algorithms for different matrix size. \cite{lrcc} proves that a sufficient condition to complete an $n\times n$ rank-$r$ matrix is $d=\mathcal{O} (n^{1.2}r\log_{10}n)$ observations. Therefore, we set $r=10$, and test $m=n=[100:200:1900]$ with $d=10 \left \lceil n^{1.2}r\log_{10}n \right \rceil$. 

\begin{table}[hbt!]
    \centering
   \resizebox{\textwidth}{!}{
    \begin{tabular}{c|ccc|ccc|ccc}
 \hline
\multirow{2}*{$n$} & \multicolumn{3}{c|}{$\text{SNR}_{\text{r}}$} & \multicolumn{3}{c|}{Iter} & \multicolumn{3}{c}{Time (seconds)} \\
           & NIHT & NN-ADMM & RC-ADMM & NIHT & NN-ADMM & RC-ADMM & NIHT & NN-ADMM & RC-ADMM \\ \hline
           100 & 27.10 & 27.10 & 27.10 & 5 & 50 & 11 & 0.01 & 0.07 & 0.01 \\
           300 & 31.46 & 31.46 & 31.46 & 5 & 39 & 12 & 0.02 & 0.52 & 0.03 \\
           500 & 33.22 & 33.22 & 33.22 & 5 & 35 & 13 & 0.04 & 1.21 & 0.08 \\
           700 & 34.33 & 34.33 & 34.33 & 5 & 32 & 14 & 0.08 & 2.32 & 0.16 \\
           900 & 35.10 & 35.10 & 35.10 & 6 & 30 & 14 & 0.14 & 3.86 & 0.21 \\
           1100 & 35.67 & 35.67 & 35.67 & 6 & 29 & 15 & 0.22 & 6.13 & 0.39 \\
           1300 & 36.13 & 36.13 & 36.13 & 6 & 28 & 16 & 0.31 & 9.27 & 0.64 \\
           1500 & 36.51 & 36.51 & 36.51 & 6 & 27 & 16 & 0.42 & 13.63 & 0.85 \\
           1700 & 36.82 & 36.82 & 36.82 & 6 & 27 & 17 & 0.59 & 19.39 & 1.40 \\
           1900 & 37.07 & 37.07 & 37.07 & 6 & 28 & 17 & 0.85 & 28.15 & 2.14 \\
    \hline
    \end{tabular}
    }
    \caption{Comparison of the algorithms with different $n$: $\text{SNR}_{\text{m}}=20$, $m=n$, $r=10$, $d=10 \left \lceil n^{1.2}r\log_{10}n \right \rceil$}
    \label{tab:4}
\end{table}

Table \ref{tab:4} shows that all the algorithms perform good though the noise exists, and NIHT converges really fast in this low-rank case. The running time of NN-ADMM increases rapidly with the matrix size increasing, which means it is not suitable for the high-dimension matrix completion. Though RC-ADMM needs more time to recover the matrix than NIHT in this experiment, when considering the case of lower sampling rate or higher rank, RC-ADMM performs better than NIHT in recovery precision. While RC-ADMM is fast enough to recover the matrix, it is also suitable for the high-dimensional matrix completion, especially for the lower sampling rate or higher rank case.

\section{Conclusions}\label{sec:8}

In this paper, we propose an ADMM algorithm for solving the general rank-constrained matrix sensing problem in both noiseless and noisy cases. And we prove the global convergence of the proposed algorithm in theory. Specialized for rank-constrained matrix completion problem, numerical experiments show that Algorithm \ref{alg:2} has an advantage in reconstruction precision comparing with NIHT, both of which are apparently better than the NN-ADMM. Though NIHT is the most efficient among the tested algorithms, Algorithm \ref{alg:2} is also fast enough. And Algorithm \ref{alg:2} shows a clear advantage in recovering when the sampling rate is really low, which is the key in matrix completion.

There are still some work to be studied in the future. First, the measurement operator is assumed to be linear, a natural question is how to extend the algorithm to the nonlinear operator case. Second, when updating $X$, we require to solve a matrix equation (\ref{equ:4.6}), how can we use the special structure of the operator $\mathcal{A}$ to avoid solve (\ref{equ:4.10}) directly. The question is closely related to the improvement of the recovery speed. Third, when updating $Y$, it is unavoidable to calculate a truncated SVD, which might cost lots of time when the matrix is high-dimensional and high-rank. We can apply some random algorithms to compute an approximate SVD \cite{random} faster. Last, though the algorithm is not sensitive to the penalty parameter $\mu$, we can still accelerate it by applying the self-adaptive penalty technique \cite{selfad}.

\bibliographystyle{unsrt}
\bibliography{references}

\end{document}